 \documentclass[11pt]{article}

 \input{amssymb.sty}

 \newtheorem{lemma}{Lemma}[section]
 \newtheorem{theorem}{Theorem}[section]
 \newtheorem{corollary}{Corollary}[section]

 \def\blemma{\begin{lemma}\sl{}\def\elemma{\end{lemma}}}
 \def\btheorem{\begin{theorem}\sl{}\def\etheorem{\end{theorem}}}

 \def\beqlb{\begin{eqnarray}}\def\eeqlb{\end{eqnarray}}
 \def\beqnn{\begin{eqnarray*}}\def\eeqnn{\end{eqnarray*}}

 \def\qed{\hfill\mbox{$\square$}\medskip}

 \def\<{\langle}\def\>{\rangle}

 \def\supp{{\mbox{\rm supp}}}

 \def\E{\mbox{\boldmath $E$}}
 \def\P{\mbox{\boldmath $P$}}
 \def\Q{\mbox{\boldmath $Q$}}
 \def\W{W}

 \def\B{{\cal B}}
 \def\J{{\cal J}}\def\L{{\cal L}}

 \def\itOmega{{\it\Omega}}\def\itPhi{{\it\Phi}}
 \def\itPsi{{\it\Psi}}

 \def\IR{\mathbb{R}}

 \def\lin{------------------}

\begin{document}

\noindent{Published in: {\it Probability Theory and Related
Fields} {\bf 127} (2003), 37--61.}

\medskip

\centerline{\bf\lin\lin\lin\lin\lin\lin}

\bigskip

\noindent{\LARGE\bf Construction of Immigration Superprocesses}

\medskip
\noindent{\LARGE\bf with Dependent Spatial Motion}

\medskip
\noindent{\LARGE\bf from One-Dimensional Excursions}

\bigskip
\noindent{Donald A. Dawson\footnote{Supported by an NSERC Research
Grant and a Max Planck Award.}}

\smallskip
\noindent{School of Mathematics and Statistics, Carleton University,}

\noindent{1125 Colonel By Drive, Ottawa, Canada K1S 5B6}

\noindent{e-mail: \tt ddawson@math.carleton.ca}

\bigskip
\noindent{Zenghu Li\,\footnote{Supported by the NSFC
(No.~10121101 and No.~10131040).}}

\smallskip
\noindent{Department of Mathematics, Beijing Normal University,}

\noindent{Beijing 100875, People's Republic of China}

\noindent{e-mail: \tt lizh@email.bnu.edu.cn}

\bigskip\bigskip

{\bf Abstract.} A superprocess with dependent spatial motion and
interactive immigration is constructed as the pathwise unique
solution of a stochastic integral equation carried by a stochastic
flow and driven by Poisson processes of one-dimensional
excursions.

\bigskip


\bigskip


\bigskip

\centerline{\bf\lin\lin\lin\lin\lin\lin}

\bigskip


\section{Introduction}

\setcounter{equation}{0}

Let $M(\IR)$ denote the space of finite Borel measures on $\IR$
endowed with a metric compatible with its topology of weak
convergence. Let $C(\IR)$ be the set of bounded continuous
functions on $\IR$. For $f\in C(\IR)$ and $\mu\in M(\IR)$ set
$\<f,\mu\> = \int fd\mu$. Let $\{w_t: t\ge0\}$ denote the
coordinate process of $C([0,\infty), M(\IR))$, which is furnished
with the locally uniform convergence. Suppose that $h$ is a
continuously differentiable function on $\IR$ such that both
$h$ and $h^\prime$ are square-integrable. Then the function
 \beqlb\label{1.1}
\rho(x) = \int_{\IR}h(y-x)h(y) dy,
\quad x\in\IR,
 \eeqlb
is twice continuously differentiable with bounded derivatives
$\rho^\prime$ and $\rho^{\prime\prime}$. Let $\sigma>0$ be a
constant. Based on the results of Dawson {\it et al} \cite{DLW01} and
Wang \cite{W98}, we shall prove that for each $\mu\in M(\IR)$
there is a unique Borel probability measure $\Q_\mu$ on
$C([0,\infty), M(\IR))$ such that, for each $\phi\in C^2(\IR)$,
 \beqlb\label{1.2}
M_t(\phi) = \<\phi,w_t\> - \<\phi,\mu\>
- \frac{1}{\,2\,}\rho(0) \int_0^t \<\phi^{\prime\prime},w_s\> ds,
\quad t\ge0,
 \eeqlb
under $\Q_\mu$ is a continuous martingale with quadratic variation
process
 \beqlb\label{1.3}
\<M(\phi)\>_t = \int_0^t\<\sigma\phi^2,w_s\> ds
+ \int_0^t ds\int_{\IR} \<h(z - \cdot) \phi^\prime, w_s\>^2 dz.
 \eeqlb
The system $\{\Q_\mu: \mu\in M(\IR)\}$ defines a diffusion
process, which we shall call a {\it superprocess with dependent
spatial motion} (SDSM). Here $\rho(0)$ is the migration rate and
$\sigma$ is the branching rate. The only difference
between the SDSM and the super Brownian motion in $M(\IR)$ is the
second term on the right hand side of (\ref{1.3}), which comes
from the dependence of the spatial motion. Because of the
dependent spatial motion, the SDSM has properties rather different
from those of the superprocess with independent spatial motion. It
is well-known that the super Brownian motion started with an
arbitrary initial state enters immediately the space of absolutely
continuous measures and its density process satisfies a stochastic
differential equation; see Konno and Shiga \cite{KS88} and Reimers
\cite{R89}. On the contrary, the SDSM lives in the space of purely
atomic measures; see Wang \cite{W97} and Theorem~\ref{t3.4} of this
paper.

The main purpose of this paper is to construct a class of
immigration diffusion processes associated with the SDSM. Let $m$
be a non-trivial $\sigma$-finite Borel measure on $\IR$ and let
$q$ be a Borel function on $M(\IR)\times \IR$ satisfying certain
regularity conditions to be specified. A modification of the SDSM
is to replace (\ref{1.2}) by
 \beqlb\label{1.4}
M_t(\phi) = \<\phi,w_t\> - \<\phi,\mu\>
- \frac{1}{\,2\,}\rho(0)\int_0^t \<\phi^{\prime\prime},w_s\> ds
- \int_0^t \<\phi q(w_s,\cdot),m\> ds,
\quad t\ge0.
 \eeqlb
A solution of the martingale problem given by (\ref{1.3}) and
(\ref{1.4}) can be interpreted as an SDSM with interactive
immigration determined by $q(w_s,\cdot)$ and the reference measure
$m$. Because of the dependence of $q(w_s,\cdot)$ on $w_s$, the
duality method of Dawson {\it et al} \cite{DLW01} and Wang \cite{W98}
fails and the uniqueness of the solution becomes a difficult
problem. The next paragraph describes our approach to the
construction of the immigration SDSM as a {\it diffusion process}.

Let $\W = C([0,\infty), \IR^+)$ and let $\tau_0(w) = \inf\{s>0:
w(s)=0\}$ for $w\in \W$. Let $\W_0$ be the set of paths $w \in \W$
such that $w(0)=w(t)=0$ for $t\ge\tau_0(w)$. We endow $\W$ and
$\W_0$ with the topology of locally uniform convergence. Let
$\Q_\kappa$ be the excursion law of the Feller branching diffusion
defined by (\ref{2.7}). Let $W(dt,dy)$ be a time-space white noise
on $[0,\infty) \times \IR$ based on the Lebesgue measure; see e.g.
Walsh \cite{W86}. Let $N_0(da,dw)$ be a Poisson random measure on $\IR
\times \W_0$ with intensity $\mu(da) \Q_\kappa(dw)$ and
$N(ds,da,du,dw)$ a Poisson random measure on $[0,\infty) \times
\IR \times [0,\infty) \times \W_0$ with intensity $ds m(da)du
\Q_\kappa(dw)$. We assume that $\{W(dt,dy)\}$, $\{N_0(da,dw)\}$
and $\{N(ds,da,du,dw)\}$ are defined on a complete standard
probability space and are independent of each other. By Dawson
{\it et al} \cite[Lemma~3.1]{DLW01} or Wang \cite[Lemma~1.3]{W97},
for any $r\ge0$ and $a\in \IR$ the stochastic equation
 \beqlb\label{1.5}
x(t) = a + \int_r^t\int_{\IR} h(y-x(s))W(ds,dy), \quad t\ge r,
 \eeqlb
has a unique continuous solution $\{x(r,a,t): t\ge r\}$, which is
a Brownian motion with quadratic variation $\rho(0)dt$. Clearly,
the system $\{x(r,a,t): t\ge r; a\in\IR\}$ determines an isotropic
stochastic flow. Let us consider the following equation:
 \beqlb\label{1.6}
Y_t
&=& \int_{\IR}\int_{\W_0} w(t) \delta_{x(0,a,t)}
N_0(da,dw)  \nonumber \\
& & + \int_0^t\int_{\IR}\int_0^{q(Y_s,a)}\int_{\W_0} w(t-s)
\delta_{x(s,a,t)}N(ds,da,du,dw),
\quad t> 0.
 \eeqlb
Our main result is that the above stochastic equation has a
pathwise unique continuous solution $\{Y_t: t>0\}$ and, if we set
$Y_0 = \mu$, then $\{Y_t: t\ge0\}$ is a diffusion process. We
prove that the distribution of $\{Y_t: t\ge0\}$ on $C([0,\infty),
M(\IR))$ solves the martingale problem (\ref{1.3}) and
(\ref{1.4}). The application of equation (\ref{1.6}) is essential
in the construction of the immigration diffusion process since the
uniqueness of solution of the martingale problem given by
(\ref{1.3}) and (\ref{1.4}) still remains open. Our equation
(\ref{1.6}) also provides useful information on the structures of
the sample paths of the immigration diffusion. For instance, from
this formulation we can immediately read of the following
properties: \begin{itemize}

\item[(i)] for any initial state $\mu\in M(\IR)$, the process
$\{Y_t: t>0\}$ consists of at most countably many atoms;

\item[(ii)] the spatial motion of the atoms of $\{Y_t: t \ge 0\}$ are
determined by the flow $\{x(r,a,t): t\ge r; a\in\IR\}$;

\item[(iii)] the mass changes of the atoms are described by
excursions of the Feller branching diffusion;

\item[(iv)] the immigration times and locations as well as the mass
excursions are selected by the Poisson random measures $N_0(da,dw)$
and $N(ds,da,du,dw)$;

\item[(v)] there are infinitely but countably many immigration
times in any non-trivial time interval if $\<q(\nu,\cdot),m\>
>0$ for all $\nu\in M(\IR)$;

\item[(vi)] for any constants $t>r>0$ there are at most only a finite
number of atoms which have lived longer than $r$ before time $t$.
\end{itemize}

A class of one-dimensional immigration diffusions was constructed
in Pitman and Yor \cite{PY82} as sums of excursions selected by Poisson
point processes. Similar constructions in infinite-dimensional
setting were discussed in Fu and Li \cite{FL03}, Li \cite{L02}, Li and
Shiga \cite{LS95} and Shiga \cite{S90}. In particular, under stronger
conditions on $q(\cdot,\cdot)$, Shiga \cite[Corollary~5.3]{S90}
constructed purely atomic measure-valued immigration branching
diffusions {\it without spatial motion} as the unique solution of
equation (\ref{1.6}) with $\delta _{x(0,a,t)}$ and $\delta
_{x(s,a,t)}$ replaced by $\delta_a$. An extension of his result to
non-trivial {\it independent spatial motion} was given recently in
Fu and Li \cite{FL03} by considering measure-valued excursions. In the
present situation, we have to put up with the hardship brought
about by the {\it dependent spatial motion}. Our approach consists
of several parts. In Section~2, we prove some useful
characterizations of the excursions of Feller branching diffusion
processes. In Section~3, we treat the case where $\<1,\mu\>
>0$ and $\<1,m\> =0$, i.e., we construct the SDSM without
immigration by excursions. The results provide useful insights
into the sample path structures of the SDSM and serve as
preliminaries of the construction of immigration processes. In
particular, from our construction it follows immediately that the
SDSM is purely atomic and, under suitable conditional
probabilities, the initial positions of its atoms are i.i.d.\
random variables with distribution $\<1,\mu\>^{-1} \mu(dx)$. In
Section~4, we consider the case where $0< \<1,m\> <\infty$ and
$q(\cdot,\cdot) \equiv 1$, i.e., the case of deterministic
immigration rate. In this case, the right hand side of (\ref{1.6})
is actually independent of $Y_s$. The problem is to show the
process $\{Y_t: t\ge0\}$ defined by this formula is a diffusion
process and solves the martingale problem (\ref{1.4}) and
(\ref{1.3}) with $q(\cdot,\cdot) \equiv 1$. This is not so easy
since the immigrants have dependent spatial motion and come
infinitely many times at any non-trivial time interval. Because of
the dependent spatial motion we cannot compute the Laplace
functional of $\{Y_t: t\ge0\}$ as in \cite{FL03} and
\cite{S90}. To resolve the difficulty, we chop off from every
excursion a part with length $1/n$ and construct a right
continuous strong Markov process $\{X^{(n)}_t: t\ge0\}$, which is
characterized as a SDSM with positive jumps. Then we obtain the
desired martingale characterization of $\{Y_t: t\ge0\}$ from that
of $\{X^{(n)}_t: t\ge0\}$ by letting $n\to\infty$. Based on those
results, we construct a pathwise unique solution of the general
equation (\ref{1.6}) in Section~5 using the techniques developed
in \cite{FL03} and \cite{S90}.


\section{Excursions of Feller branching diffusions}

\setcounter{equation}{0}

Let $\beta>0$ be a constant and $\{B(t): t\ge0\}$ a standard
Brownian motion. For any initial condition $\xi(0) = x \ge0$ the
stochastic differential equation
 \beqlb\label{2.1}
d \xi(t) = \sqrt{\beta\xi(t)} d B(t),
\quad t\ge0,
 \eeqlb
has a unique solution $\{\xi(t): t\ge0\}$, which is a diffusion
process on $[0,\infty)$. The transition semigroup $(Q_t)_{t\ge0}$
of the process is determined by
 \beqlb\label{2.2}
\int_0^\infty e^{-z y}Q_t(x,dy)
=
\exp\{- xz(1+\beta tz/2)^{-1}\},
\quad t,x,z\ge0;
 \eeqlb
see e.g. Ikeda and Watanabe \cite[p.236]{IW89}. In this paper, we call
any diffusion process $\{\xi(t): t\ge0\}$ a {\it Feller branching
diffusion with constant branching rate $\beta$}, or simply a {\it
$\beta$-branching diffusion} if it has transition semigroup
$(Q_t)_{t\ge0}$. Letting $z \to \infty$ in (\ref{2.2}) we get
 \beqlb\label{2.3}
Q_t (x,\{0\})
=
\exp\{-2x/\beta t\},
\quad t>0,x\ge0.
 \eeqlb
In view of (\ref{2.1}), $\{\xi(t): t\ge0\}$ is a continuous
martingale with quadratic variation $\beta \xi(t)dt$. In general,
if $\{\eta(t): t\ge0\}$ is a continuous martingale with quadratic
variation $\sigma(t) \eta(t)dt$ for a predictable process
$\{\sigma(t): t\ge0\}$, we call it a {\it Feller branching
diffusion with branching rate} $\{\sigma(t): t\ge0\}$.

Let $Q_t^\circ(x,\cdot)$ denote the restriction of the measure
$Q_t (x,\cdot)$ to $(0,\infty)$. Since the origin $0$ is a trap
for the $\beta$-branching diffusion process, the family of kernels
$(Q_t^\circ) _{t\ge0}$ also constitute a semigroup. In view of the
infinite divisibility implied by (\ref{2.2}), there is a family of
{\it canonical measures} $(\kappa_t)_{t>0}$ on $(0,\infty)$ such
that
 \beqlb\label{2.4}
\int_0^\infty (1- e^{-zy})\kappa_t(dy)
=
z(1+\beta tz/2)^{-1},
\quad t>0,z\ge0.
 \eeqlb
Indeed, we have
 \beqlb\label{2.5}
\kappa_t(dy)
= 4(\beta t)^{-2}e^{-2y/\beta t}dy,
\quad t>0, x>0.
 \eeqlb
Based on (\ref{2.2}) and (\ref{2.4}) one may check that
 \beqlb\label{2.6}
\int_0^\infty (1- e^{-zy})\kappa_{r+t}(dy)
=
\int_0^\infty \kappa_r(dy) \int_0^\infty (1- e^{-zy})
Q_t^\circ (x,dy),
\quad r,t>0,z\ge0.
 \eeqlb
Then $\kappa_rQ_t^\circ = \kappa_{r+t}$ and hence $(\kappa_t)
_{t>0}$ is an {\it entrance law} for $(Q_t^\circ) _{t\ge0}$. It is
known that there is a unique $\sigma$-finite measure $\Q_\kappa$
on $(\W_0,\B(\W_0))$ such that
 \beqlb\label{2.7}
\Q_\kappa\{w(t_1)\in dy_1,\cdots,w(t_n)\in dy_n\}
=
\kappa_{t_1}(dy_1)Q^\circ_{t_2-t_1}(y_1,dy_2)\cdots
Q^\circ_{t_n-t_{n-1}}(y_{n-1},dy_n)
 \eeqlb
for $0<t_1<t_2<\cdots<t_n$ and $y_1,y_2,\cdots,y_n\in (0,\infty)$;
see e.g.\ Pitman and Yor \cite{PY82} for details. The measure
$\Q_\kappa$ is known as the {\it excursion law} of the
$\beta$-branching diffusion. Let ${\cal B}_t = {\cal B}_t(\W_0)$
denote the $\sigma$-algebra on $\W_0$ generated by $\{w(s): 0\le
s\le t\}$. Roughly speaking, (\ref{2.7}) asserts that $\{w(t):
t>0\}$ is a $\beta$-branching diffusion relative to
$(\Q_\kappa,{\cal B}_t)$ with one-dimensional distributions
$\{\kappa_t: t>0\}$.

For $r>0$, let $\tau_r(w) = r\vee\tau_0(w)$, and let
$\Q_{\kappa,r}$ denote the restriction of $\Q_\kappa$ to $\W_r :=
\{w\in \W_0: \tau_0(w)>r\}$. Observe that
 \beqlb\label{2.8}
\Q_\kappa(\W_r)
=
\Q_{\kappa,r}(\W_r)
= \kappa_r(0,\infty)
= 2/\beta r,
\quad r>0.
 \eeqlb
The following theorem gives a stochastic equation for the
excursions.

\btheorem\label{t2.1} For any $r>0$, the coordinate process
$\{w(t): t\ge r\}$ under $\Q_{\kappa,r} \{\,\cdot\,| {\cal B}_r\}$
is a $\beta$-branching diffusion. Moreover, there is a measurable
mapping $B: \W_0 \to \W$ such that $\{B(w,t)\equiv B(w)(t): t\ge
r\}$ under $\Q_{\kappa,r} \{\,\cdot\,| {\cal B}_r\}$ is a Brownian
motion stopped at time $\tau_r(w)$ and
 \beqlb\label{2.9}
dw(t)
=
\sqrt{\beta w(t)} dB(w,t),
\quad t\ge r.
 \eeqlb
\etheorem

\noindent{\it Proof.} By (\ref{2.7}), $\{w(t): t\ge r\}$ under
$\Q_{\kappa,r}\{\,\cdot\,|{\cal B}_r\}$ is a $\beta$-branching
diffusion or, equivalently, a continuous martingale with quadratic
variation $\beta w(t)dt$. Thus
 \beqlb\label{2.10}
\tilde B_r(w,t)
:=
\int_r^{t\land \tau_r}\frac{1}{\sqrt{\beta w(s)}} dw(s),
\quad t\ge r,
 \eeqlb
under $\Q_{\kappa,r} \{\,\cdot\,|{\cal B}_r\}$ is a Brownian
motion stopped at time $\tau_r(w)$. We may define another Markov
transition semigroup $(Q_t^h)_{t\ge0}$ on $[0,\infty)$ by
 \beqnn
Q_t^h(x,dy)
=
\left\{\begin{array}{ll}
x^{-1}yQ_t(x,dy)  \quad &\mbox{for $t>0, x>0, y\ge 0$,}  \\
y\kappa_t(dy)                  &\mbox{for $t>0, x=0, y\ge 0$,}
\end{array}\right.
 \eeqnn
which is a Doob's $h$-transform of $(Q_t)_{t\ge0}$. It is not hard
to check that $\bar{\Q}_{\kappa,r}(dw) := w(r) \Q_{\kappa,r} (dw)$
defines a probability measure. Moreover, (\ref{2.7}) implies that
$\{w(t): 0\le t\le r\}$ under $\bar{\Q}_{\kappa,r}$ is a diffusion
process with transition semigroup $(Q_t^h)_{t\ge0}$ and generator
$2^{-1}\beta xd^2/dx^2 + \beta d/dx$. Thus $\{2 \sqrt{w(t)/\beta}:
0\le t\le r\}$ under $\bar{\Q}_{\kappa,r}$ is a $4$-dimensional
Bessel diffusion. It follows that
 \beqlb\label{2.11}
m(t) := w(t) - \beta t,
\quad 0\le t\le r,
 \eeqlb
is a continuous martingale with quadratic variation $\beta
w(t)dt$ and the limits
 \beqlb\label{2.12}
\bar M_r(w,t)
:=
\lim_{u\to0^+}\int_u^t \frac{1}{\sqrt{\beta w(s)}} dm(s),
\quad 0<t\le r,
 \eeqlb
exist in $L^2(\bar{\Q}_{\kappa,r})$-sense. Let $\bar M_r(w,0)=0$.
Then $\{\bar M_r(w,t): 0\le t\le r\}$ under $\bar{\Q}_{\kappa,r}$ is
a Brownian motion. By Shiga and Watanabe \cite[Theorem~3.3.ii]{SW73},
for any $0<\epsilon<1$ we have $\bar{\Q}_{\kappa,r} \{\sqrt{w(t)} <
t^{(1+\epsilon)/2}$ infinitely often as $t\to0^+\} =0$. By
(\ref{2.11}) and (\ref{2.12}), the limits
 \beqnn
\bar B_r(w,t)
:=
\lim_{u\to0^+}\int_u^t \frac{1}{\sqrt{\beta w(s)}} dw(s),
\quad 0<t\le r,
 \eeqnn
also exist in $L^2(\bar{\Q}_{\kappa,r})$-sense. Indeed, setting $\bar
B_r(w,0)=0$ we have
 \beqnn
\bar B_r(w,t)
=
\bar M_r(w,t) + \int_{0^+}^t \frac{\sqrt{\beta}}{\sqrt{w(s)}} ds,
\quad 0\le t\le r.
 \eeqnn
Now let $B(w,t) = \bar B_r(w,t)$ for $0\le t\le r$ and $B(w,t) =
\bar B_r(w,r) + \tilde B_r(w,t)$ for $t\ge r$. Since
$\Q_{\kappa,r}$ is absolutely continuous relative to $\bar{\Q}
_{\kappa,r}$, we see that $B(w,t)$ is uniquely defined on $\W_r$
out of a $\Q_\kappa$-null set and it does not depend on the
particular choice of $r>0$. Using this property, we can extend the
definition of $B(w,t)$ to the whole space $\W_0$ so that $\{B(w,t):
t\ge r\}$ under $\Q_{\kappa,r} \{\,\cdot\,|{\cal B}_r\}$ is a
Brownian motion stopped at time $\tau_r(w)$ and satisfies
(\ref{2.9}). \qed


\section{SDSM without immigration}

\setcounter{equation}{0}

In this section, we give a rigorous construction of a purely
atomic version of the SDSM with an arbitrary initial state. The
results are useful in our study of the associated immigration
processes. We consider a general branching density $\sigma(\cdot)
\in C(\IR)^+$ and assume there is a constant $\epsilon>0$ such
that $\sigma (x)\ge \epsilon$ for all $x\in \IR$. From the results
in Dawson {\it et al} \cite{DLW01} and Wang \cite{W98} we know that
the generator $\L$ of the SDSM with branching density $\sigma
(\cdot)$ is expressed as
 \beqlb\label{3.1}
\L F(\nu)
&=& \frac{1}{2}\rho(0)\int_{\IR}\frac{d^2}{dx^2}
\frac{\delta F(\nu)}{\delta\nu(x)}\nu(dx)  \nonumber  \\
& &
+\,\frac{1}{2}\int_{\IR^2}\rho(x-y)
\frac{d^2}{dxdy}\frac{\delta^2 F(\nu)}
{\delta\nu(x)\delta\nu(y)}\nu(dx)\nu(dy)  \nonumber \\
& &
+\,\frac{1}{2} \int_{\IR}\sigma(x) \frac{\delta^2 F(\nu)}
{\delta\nu(x)^2}\nu(dx).
 \eeqlb
The domain of $\L$ include functions on $M(\IR)$ of the form
$F_{n,f}(\nu) := \int fd\nu^n$ with $f\in C^2(\IR^n)$ and
functions of the form
 \beqlb\label{3.2}
F_{f,\{\phi_i\}}(\nu) := f(\<\phi_1,\nu\>, \cdots, \<\phi_n,\nu\>)
 \eeqlb
with $f\in C^2(\IR^n)$ and $\{\phi_i\}\subset C^2(\IR)$. Let
${\cal D}(\L)$ denote the collection of all those functions.

\btheorem\label{t3.1} {\rm (Dawson {\it et al}, 2001; Wang, 1998)} For
each $\mu\in M(\IR)$ there is a unique probability measure
$\Q_\mu$ on the space $C([0,\infty), M(\IR))$ such that $\Q_\mu
\{w_0 = \mu\} =1$ and under $\Q_\mu$ the coordinate process
$\{w_t: t\ge0\}$ solves the $(\L$, ${\cal D}(\L))$-martingale
problem. Therefore, the system $\{\Q_\mu: \mu\in M(\IR)\}$ defines
a diffusion process generated by the closure of $(\L$, ${\cal
D}(\L))$. \etheorem

\btheorem\label{t3.2} A probability measure $\Q_\mu$ on
$C([0,\infty), M(\IR))$ is a solution of the $(\L$, ${\cal
D}(\L))$-martingale problem with $\Q_\mu \{w_0=\mu\} =1$ if and
only if for each $\phi\in C^2(\IR)$,
 \beqlb\label{3.3}
M_t(\phi) = \<\phi,w_t\> - \<\phi,\mu\>
- \frac{1}{\,2\,}\rho(0) \int_0^t \<\phi^{\prime\prime},w_s\> ds,
\quad t\ge0,
 \eeqlb
under $\Q_\mu$ is a continuous martingale with quadratic variation
process
 \beqlb\label{3.4}
\<M(\phi)\>_t = \int_0^t\<\sigma\phi^2,w_s\> ds
+ \int_0^t ds\int_{\IR} \<h(z - \cdot) \phi^\prime, w_s\>^2 dz.
 \eeqlb
\etheorem

\noindent{\it Proof.} Suppose that $\Q_\mu$ is a probability
measure on $C([0,\infty), M(\IR))$ such that $\Q_\mu \{w_0 =
\mu\} =1$ and
 \beqlb\label{3.5}
F(w_t) - F(w_0) - \int_0^t \L F(w_s) ds,
\quad t\ge0,
 \eeqlb
is a continuous martingale for every $F\in{\cal D}(\L)$. Comparing
the martingales related to the functions $\mu \mapsto
\<\phi,\mu\>$ and $\mu\mapsto \<\phi,\mu\>^2$ and using It\^o's
formula we see that (\ref{3.3}) is a continuous martingale with
quadratic variation process (\ref{3.4}). Conversely, suppose that
$\Q_\mu$ is a probability measure on $C([0,\infty), M(\IR))$ under
which (\ref{3.3}) is a continuous martingale with quadratic
variation process (\ref{3.4}) for each $\phi\in C^2(\IR)$. Observe
that for the function $F_{f,\{\phi_i\}}$ defined by (\ref{3.2}) we
have
 \beqnn
\L F_{f,\{\phi_i\}}(\nu)
&=&
\frac{1}{2}\rho(0) \sum_{i=1}^nf_i^\prime(\<\phi_1,\nu\>,\cdots,
\<\phi_n,\nu\>) \<\phi_i^{\prime\prime},\nu\>  \nonumber  \\
& &
+ \frac{1}{2}\sum_{i,j=1}^n f_{ij}^{\prime\prime}(\<\phi_1,\nu\>,
\cdots, \<\phi_n,\nu\>) \int_{\IR^2}\rho(x-y)\phi_i^\prime(x)
\phi_j^\prime(y) \nu(dx)\nu(dy)   \nonumber \\
& &
+ \frac{1}{2}\sum_{i,j=1}^n f_{ij}^{\prime\prime} (\<\phi_1,\nu\>,
\cdots, \<\phi_n,\nu\>) \<\sigma\phi_i\phi_j,\nu\>.
 \eeqnn
By It\^o's formula we see that (\ref{3.5}) is a continuous
martingale if $F=F_{f,\{\phi_i\}}$. Then the theorem follows by an
approximation of an arbitrary $F\in{\cal D}(\L)$. \qed

We now consider the construction of the trajectories of the SDSM.
Suppose that $(\itOmega,{\cal F},\P)$ is a complete standard
probability space on which we have a white noise $\{W(ds,dy)\}$ on
$[0,\infty) \times \IR$ based on the Lebesgue measure; see e.g.\
Walsh \cite{W86}. By Dawson {\it et al} \cite[Lemma~3.1]{DLW01} or
Wang \cite[Lemma~1.3]{W97}, for any $a\in \IR$ the equation
 \beqlb\label{3.6}
x(t) = a + \int_0^t\int_{\IR} h(y-x(s))W(ds,dy),
\quad t\ge 0,
 \eeqlb
has a unique solution $\{x(a,t): t\ge 0\}$. For a fixed constant
$\beta>0$ let
 \beqlb\label{3.7}
\psi(a,t) = \beta^{-1}\int_0^t \sigma(x(a,s)) ds,
\quad t\ge 0, a\in\IR.
 \eeqlb
Suppose we also have on $(\itOmega,{\cal F},\P)$ a sequence of
independent $\beta$-branching diffusions $\{\xi_i(t): t\ge0;
i=1,2,\cdots\}$ independent of $\{W(ds,dy)\}$. We assume each
$\xi_i(0) \ge0$ is deterministic and $\sum_{i=1}^\infty \xi_i(0) <
\infty$. Let $\{a_i: i=1,2,\cdots\}$ be a sequence of real
numbers. Let $\xi_i(a,t) = \xi_i (\psi (a,t))$ and let ${\cal
G}_t$ be the $\sigma$-algebra generated by all $\P$-null
sets and the family of random variables
 \beqlb\label{3.8}
\{W([0,s]\times B): 0\le s\le t; B\in{\cal B}(\IR)\}
\quad\mbox{and}\quad
\{\xi_i(a_i,s): 0\le s\le t; i=1,2,\cdots\}.
 \eeqlb

For $r>0$ let $n(r) = \#\{i: \xi_i(r)>0\}$ and $n^\psi(r) = \#
\{i: \xi_i(a_i,r)>0\}$, where $\#\{\cdots\}$ denotes the number of
elements of the set $\{\cdots\}$. Since zero is a trap for the
$\beta$-branching diffusion, both $n(r)$ and $n^\psi(r)$ are a.s.\
non-increasing in $r>0$.

\blemma\label{l3.1} We have $n(r) < \infty$ and $n^\psi(r) <
\infty$ a.s.\ for each $r>0$. \elemma

\noindent{\it Proof.} In view of (\ref{2.3}), we have
 \beqnn
\sum_{i=1}^\infty \P\{\xi_i(r)>0\}
=
\sum_{i=1}^\infty [1-\exp\{-2\xi_i(0)/\beta r\}]
\le
\frac{2}{\beta r}\sum_{i=1}^\infty \xi_i(0)
< \infty.
 \eeqnn
Then an application of the Borel-Cantelli lemma yields that a.s.\
$n(r) <\infty$. By (\ref{3.7}) and the assumption $\sigma(x)\ge
\epsilon$, we have $\psi (a_i,t) \ge \epsilon t/\beta$. But
$0$ is a trap for $\{\xi_i(t): t\ge0\}$, so $n^\psi (r) <
n(\epsilon r/\beta) < \infty$ a.s.\ for each $r>0$. \qed

\btheorem\label{t3.3} The process $\{X_t: t\ge0\}$ defined by
 \beqlb\label{3.9}
X_t = \sum_{i=1}^\infty \xi_i(a_i,t)\delta_{x(a_i,t)},
\quad t\ge0,
 \eeqlb
relative to $({\cal G}_t)_{t\ge0}$ is an SDSM.
\etheorem

\noindent{\it Proof.} By the assumption of independence,
$\{\xi_i(t): t\ge0; i=1,2,\cdots\}$ and $\{W([0,t] \times B):
t\ge0, B \in {\cal B}(\IR)\}$ are martingales relative to the
filtration $({\cal G}_t)_{t\ge0}$. Given $\{W(ds,dy)\}$, the
processes $\{\psi (a_i,\cdot): i=1,2,\cdots\}$ are deterministic.
Then the time-changed processes $\{\xi_i (a_i,\cdot): i= 1,2,
\cdots\}$ are independent martingales under $\P\{\cdot\,| W\}$.
Moreover, we have a.s.\
 \beqnn
\<\xi_i(a_i)\> (t)
=
\int_0^{\psi(a_i,t)}\beta\xi_i(s) ds
=
\int_0^t \beta\xi_i(a_i,u) d\psi(a_i,u)
=
\int_0^{t} \sigma(x_i(u))\xi_i(a_i,u) du
 \eeqnn
first under $\P\{\cdot\,| W\}$ and then under the non-conditional
probability $\P$. By the same reasoning we get $\<\xi_i(a_i),
\xi_j(a_j)\>(t) \equiv 0$ a.s.\ under $\P$ for $i\neq j$. By
It\^o's formula,
 \beqnn
\xi_i(a_i,t)\phi(x(a_i,t))
&=&
\xi_i(0)\phi(a_i)
+ \int_0^t\int_{\IR}\xi_i(a_i,s)\phi^\prime(x(a_i,s))
h(y-x(a_i,s))W(ds, dy)      \\
& & + \,\frac{1}{\,2\,} \int_0^tds\int_{\IR}\xi_i(a_i,s)
\phi^{\prime\prime} (x(a_i,s)) h(y-x(a_i,s))^2 dy  \\
& & + \,\int_0^t\phi(x(a_i,s)) d\xi_i(a_i,s)
 \eeqnn
for $\phi\in C^2(\IR)$. Taking the summation $\sum_{i=1}^\infty$
we get
 \beqnn
\<\phi,X_t\>
= \<\phi,X_0\> + M_t(\phi)
+ \frac{1}{\,2\,}\rho(0) \int_0^t \<\phi^{\prime\prime},X_s\> ds,
\quad t\ge0,
 \eeqnn
where
 \beqnn
M_t(\phi)
:=
\int_0^t\int_{\IR}\<h(y-\cdot)\phi^\prime,X_s\> W(ds, dy)
+ \sum_{i=1}^\infty \int_0^t\phi(x(a_i,s))d\xi_i(a_i,s),
 \eeqnn
is a continuous martingale relative to $({\cal G}_t)_{t\ge0}$ with
quadratic variation process
 \beqnn
\<M(\phi)\>_t
=
\int_0^t ds\int_{\IR} \<h(z - \cdot) \phi^\prime, X_s\>^2 dz
+ \int_0^t\<\sigma\phi^2,X_s\> ds.
 \eeqnn
Then we have the desired result by Theorem~\ref{t3.2}. \qed

By Theorem~\ref{t3.3}, if the SDSM is started from an initial
state in $M_a(\IR)$, purely atomic measures on $\IR$, it lives in
this space forever. More precisely, the position of its $i$th atom
is described by $\{x(a_i,t): t\ge0\}$ and its mass by $\{\xi_i
(a_i,t): t\ge0\}$. As observed in Wang \cite[p.756]{W98}, each
$\{x(a_i,t): t\ge0\}$ is a Brownian motion with quadratic
variation $\rho(0)dt$. If $a_i=a_j$, we have $x(a_i,t) = x(a_j,t)$
for all $t\ge0$ by the uniqueness of solution of (\ref{3.6}). On
the contrary, if $a_i\neq a_j$, then $\{x(a_i,t): t\ge0\}$ and
$\{x(a_j,t): t\ge0\}$ never hit each other.

To consider a more general initial state, we need a lemma for
Poisson random measures. Suppose that $E$ and $F$ are metrizable
topological spaces. Let $\mu \in M(E)$ and let $q(x,dy)$ be a
Borel probability kernel from $E$ to $F$. Then
 \beqlb\label{3.10}
\int_E\int_F h(x,y)\nu(dx,dy)
=
\int_E\mu(dx)\int_F h(x,y)q(x,dy),
\quad h\in C(E\times F),
 \eeqlb
defines a measure $\nu\in M(E\times F)$. Let $Y$ be a Poisson
random measure on $E\times F$ with intensity $\nu$ and let
$\eta = Y(E\times F)$. The following lemma shows that we can
recover the kernel $q(x,dy)$ from the atoms of $Y$ by a {\it
suitable enumeration}.

\blemma\label{l3.2} In the situation described above, $X(\cdot) :=
Y(\cdot\times F)$ defines a Poisson random measure on $E$ with
intensity $\mu$. Suppose that $\{(x_i,y_i): i=1,\cdots,\eta\}$ is
an enumeration of the atoms of $Y(dx,dy)$ which only uses
information from $X$. Then, given $\eta=k$ and $x_i = c_i \in E$
$(i=1,\cdots,k)$, the sequence $\{y_i: i=1,\cdots,k\}$ is formed
of independent random variables with distributions $\{q(c_i,\cdot):
i=1,\cdots,k\}$. \elemma

\noindent{\it Proof.} We first consider a special version of the
Poisson random measure $Y$ constructed as follows. Let $\eta$ be a
Poisson random variable with parameter $\mu(E)$ and $\{(x_i,y_i):
i=1,2,\cdots\}$ a sequence of i.i.d.\ random variables in $E\times
F$ which are independent of $\eta$ and have common distribution
$\mu(E)^{-1}\nu$. Then
 \beqnn
Y(dx,dy) := \sum_{i=1}^{\eta}\delta_{(x_i,y_i)}(dx,dy),
\quad x\in E, y\in F,
 \eeqnn
is a Poisson random measure $Y$ with intensity $\nu(dx,dy)$ and
 \beqnn
X(dx) := \sum_{i=1}^{\eta} \delta_{x_i}(dx),
\quad x\in E,
 \eeqnn
is a Poisson random measure on $E$ with intensity $\mu(dx)$.
Observe that, for any integer $k\ge1$, any permutation
$\{i_1,\cdots,i_k\}$ of $\{1,\cdots,k\}$ and any sequence $\{c_j:
j=1,\cdots,k\} \subset E$, under $\P\{\,\cdot\,|x_{i_j}=c_j:
j=1,\cdots,k\}$ the sequence $\{y_{i_j}: j=1,\cdots,k\}$ consists
of independent random variables with distributions
$\{q(c_j,\cdot): j= 1,\cdots,k\}$. This proves the lemma for the
special version of $Y$. The result for an arbitrary realization of
the Poisson random measure holds by the uniqueness of
distribution. \qed

Now we consider the construction of the SDSM with a general
initial state $\mu\in M(\IR)$ with $\<1,\mu\> >0$. Suppose we have
on some complete standard probability space $(\itOmega,{\cal
F},\P)$ a time-space white noise $W(ds,dy)$ on $[0,\infty) \times
\IR$ based on the Lebesgue measure and a Poisson random measure
$N(da,dw)$ on $\IR\times \W_0$ with intensity $\mu(da) \Q_\kappa
(dw)$, where $\Q_\kappa$ denotes the excursion law of the
$\beta$-branching diffusion defined by (\ref{2.7}). Assume that
$\{W(ds,dy)\}$ and $\{N(da,dw)\}$ are independent. Let $m(r) =
N(\IR \times \W_r)$ for $r>0$. In view of (\ref{2.8}), we have
 \beqnn
\E\{m(r)\}
=
\<1,\mu\>\Q_\kappa(\W_r)
=
2\<1,\mu\>/r\beta,
 \eeqnn
and hence a.s.\ $m(r) < \infty$. Thus we can enumerate the atoms
of $N(da,dw)$ into a sequence $\supp(N) =\{(a_i,w_i): i= 1,2,
\cdots\}$ such that a.s.\ $\tau_0(w_{i+1}) < \tau_0(w_i)$ for all
$i\ge1$ and $\tau_0(w_i) \to 0$ as $i\to\infty$. Clearly, $m(r)=i$
for $\tau_0(w_{i+1})\le r < \tau_0(w_i)$, and $m(r)=0$ for $r \ge
\tau_0(w_1)$. Let $\psi(a,t)$ be defined by (\ref{3.7}) and let
$w(a,t) = w(\psi(a,t))$ for $w\in \W_0$. For $r>0$ let $\supp_r(N)
= \{(a_i,w_i): i= 1, \cdots, m(r)\}$, let $\supp_r^\psi (N) =
\{(a_i,w_i): w_i(a_i,r)>0; i= 1,2, \cdots\}$ and let $m^\psi(r) =
\# \{\supp_r^\psi (N)\}$. As in the proof of Lemma~\ref{l3.2}, one
can show that a.s.\ $m^\psi(r) < \infty$. Then we have the
following

\blemma\label{l3.3} For each $r>0$, we have a.s.\ $m(r) < \infty$
and $m^\psi(r) < \infty$. \elemma

Let us see how to recover branching diffusions under some
conditional probabilities by reordering suitably the atoms of
$N(da, dw)$. For $t\ge0$ let ${\cal G}_t$ be the $\sigma$-algebra
generated by all $\P$-null sets and the families random variables
 \beqlb\label{3.11}
\{W([0,s]\times B): 0\le s\le t; B\in{\cal B}(\IR)\}
\quad\mbox{and}\quad
\{w_i(a_i,s): 0\le s\le t; i=1,2,\cdots\}.
 \eeqlb

\blemma\label{l3.4}  For each $r>0$ there is an enumeration
$\{(a_{i_j}, w_{i_j}): j= 1,\cdots, m^\psi(r)\}$ of $\supp_r^\psi
(N)$ which only uses information from ${\cal G}_r$ and satisfies
that $\{w_{i_j} (\psi(a_{i_j},r)+t): t\ge 0; j=1, \cdots,
m^\psi(r)\}$ under $\P\{\,\cdot\, |{\cal G}_r\}$ are independent
$\beta$-branching diffusions which are independent of $\{W(dt,dy):
t\ge r; y\in \IR\}$. \elemma

\noindent{\it Proof.} As observed in the proof of Lemma~\ref{l3.1}
we have $\psi(a_i,r) \ge \epsilon r/\beta$. Since the finite
measure $\kappa_{\epsilon r/\beta}(dy)$ on $(0,\infty)$ is
absolutely continuous, the elements of the random set
$\{w_i(\epsilon r/\beta): i=1,\cdots,m(\epsilon r/\beta)\}$ are
a.s.\ mutually distinct. Let ${\cal F}_t := \sigma(\{w_i(s): 0\le
s\le t; i=1,2,\cdots\})$ and let $\{(a_{k_j}, w_{k_j}):
j=1,\cdots,m(\epsilon r/\beta)\}$ be the enumeration of
$\supp_{\epsilon r/\beta}(N)$ so that $w_{k_1}(\epsilon r/\beta) <
\cdots < w_{k_{m(\epsilon r/\beta)}}(\epsilon r/\beta)$. Note that
this enumeration only uses information from ${\cal F}_{\epsilon
r/\beta}$. An application of Theorem~\ref{t2.1} and
Lemma~\ref{l3.2} shows that $\{w_{k_j}(\epsilon r/\beta+t): t\ge
0; j=1, \cdots, m(\epsilon r/\beta)\}$ under $\P\{\,\cdot\, |{\cal
F}_{\epsilon r/\beta}\}$ are independent $\beta$-branching
diffusions. By the independence of $\{W(ds,dy)\}$ and
$\{N(da,dw)\}$, we know that $\{w_{k_j}(\epsilon r/\beta+t): t\ge
0; j=1, \cdots, m(\epsilon r/\beta)\}$ are also independent
$\beta$-branching diffusions under $\P\{\,\cdot\, |{\cal
F}_{\epsilon r/\beta}^e\}$, where ${\cal F}_{\epsilon r/\beta}^e$
is the $\sigma$-algebra generated by ${\cal F}_{\epsilon r/\beta}$
and the family of random variables $\{W([0,s]\times B): 0\le s\le
r; B\in{\cal B}(\IR)\}$. Moreover, $\{w_{k_j}(\epsilon r/\beta+t):
t\ge 0; j=1, \cdots, m(\epsilon r/\beta)\}$ and $\{W(dt,dy): t\ge
r; y\in \IR\}$ are independent under $\P\{\,\cdot\, |{\cal
F}_{\epsilon r/\beta}^e\}$. Observe that ${\cal G}_r$ is generated
by ${\cal F}_{\epsilon r/\beta}^e$, all $\P$-null sets and the
family of random variables $\{w_i(s): \epsilon r/\beta\le s\le
\psi(a_i,r); i=1, \cdots, m(\epsilon r/\beta)\}$. It follows that
$\{w_{k_j}(\psi (a_{k_j},r)+t): t\ge 0; j=1, \cdots, m(\epsilon
r/\beta)\}$ under $\P\{\,\cdot\, |{\cal G}_r\}$ are also
independent $\beta$-branching diffusions which are independent of
$\{W(dt,dy): t\ge r; y\in \IR\}$. Finally, we remove the elements
of $\{(a_{k_j}, w_{k_j}): j= 1,\cdots, m(\epsilon r/\beta)\}$ with
$w (a_{k_j},r) =w_{k_j} (\psi (a_{k_j},r)) =0$ and relabel the
remaining elements to get an enumeration $\{(a_{i_j}, w_{i_j}):
j=1, \cdots, m^\psi(r)\}$ of $\supp_r^\psi (N)$ so that
$w_{i_1}(\epsilon r/\beta) < \cdots < w_{i_{m^\psi(r)}} (\epsilon
r/\beta)$. Clearly, this enumeration only uses information from
${\cal G}_r$ and has the desired property. \qed

\btheorem\label{t3.4} Let $\{X_t: t\ge0\}$ be defined by $X_0=\mu$
and
 \beqlb\label{3.12}
X_t =
\sum_{i=1}^\infty w_i(a_i,t)\delta_{x(a_i,t)}
=
\int_{\IR}\int_{\W_0} w(a,t)\delta_{x(a,t)} N(da,dw),
\quad t>0.
 \eeqlb
Then $\{X_t: t\ge0\}$ relative to $({\cal G}_t)_{t\ge0}$ is an
SDSM. \etheorem

\noindent{\it Proof.} Let $(Q_t)_{t\ge0}$ denote the transition
semigroup of the SDSM. For $r>0$ we see by Lemma~\ref{l3.4} and
Theorem~\ref{t3.3} that $\{X_t: t\ge r\}$ under $\P\{\,\cdot\,
|{\cal G}_r\}$ is a Markov process with transition semigroup
$(Q_t)_{t\ge0}$. Thus $\{X_t: t>0\}$ is a Markov process with
transition semigroup $(Q_t)_{t\ge0}$. We shall prove that the
random measure $X_t$ has distribution $Q_t(\mu,\cdot)$ for $t>0$
so that the desired result follows from the uniqueness of
distribution of the SDSM. Using the notation of the proof of
Lemma~\ref{l3.4}, we have that $\{w_{k_j}(\epsilon r/\beta+t):
t\ge 0; j=1, \cdots, m(\epsilon r/\beta)\}$ under $\P\{\,\cdot\,
|{\cal F}_{\epsilon r/\beta}\}$ are independent $\beta$-branching
diffusions which are independent of $\{W(dt,dy): t\ge 0;
y\in\IR\}$. By Theorem~\ref{t3.3},
 \beqlb\label{3.13}
X_t^{(r)}
:=
\sum_{j=1}^{m(\epsilon r/\beta)} w_{k_j}(\epsilon r/\beta
+\psi(a_{k_j},t))\delta_{x(a_{k_j},t)},
\quad t\ge0,
 \eeqlb
under $\P\{\cdot\,|{\cal F}_{\epsilon r/\beta}\}$ is an SDSM with
initial state
 \beqnn
X_0^{(r)} = \sum_{j=1}^{m(\epsilon r/\beta)} w_{k_j}
(\epsilon r/\beta)\delta_{a_{k_j}}.
 \eeqnn
This implies that $\{X_t^{(r)}: t\ge0\}$ under the non-conditional
probability $\P$ is an SDSM relative to the filtration $({\cal
H}_t^{(r)})_{t\ge0}$, where ${\cal H}_t^{(r)}$ is generated by
${\cal F} _{\epsilon r/\beta}$ and $\{X_s^{(r)}: 0\le s\le t\}$.
For any $f\in C(\IR)^+$, we have by (\ref{2.4}) that
 \beqnn
\E\exp\{-\<f,X_0^{(r)}\>\}
&=&
\E\exp\bigg\{-\int_{\IR}\int_{\W_0}w(\epsilon r/\beta)f(a)
N(da,dw)\bigg\}  \\
&=&
\exp\bigg\{-\int_{\IR}\mu(da)\int_{\W_0}(1
-e^{-w(\epsilon r/\beta)f(a)})\Q_\kappa(dw)\bigg\}  \\
&=&
\exp\bigg\{-\int_{\IR}f(a)(1+\epsilon r f(a)/2)^{-1}
\mu(da)\bigg\},
 \eeqnn
which converges to $\exp\{-\<f,\mu\>\}$ as $r\to 0$. Thus
$X_0^{(r)} \to \mu$ in distribution as $r\to 0$. Indeed, if we set
$X_0^{(0)}=\mu$, then $\{X_0^{(r)}: r\ge 0\}$ is a measure-valued
branching diffusion without migration; see
\cite[Theorem~3.6]{S90}. By the Feller property of the SDSM, the
distribution of $X_t^{(r)}$ converges to $Q_t(\mu,\cdot)$ as
$r\to0$. Since $\psi(a_{k_j},t)\ge \epsilon t/\beta$, we can
rewrite (\ref{3.13}) as
 \beqnn
X_t^{(r)}
:=
\sum_{i=1}^{m(\epsilon t/\beta)} w_i(\epsilon r/\beta
+\psi(a_i,t))\delta_{x(a_i,t)}.
 \eeqnn
Then for fixed $t>0$ we have $X_t^{(r)} \to X_t$ a.s.\ as $r\to0$
and hence $X_t$ has distribution $Q_t(\mu,\cdot)$. \qed

By Theorem~\ref{t3.4}, the SDSM started with an arbitrary initial
measure enters the space $M_a(\IR)$ of purely atomic measures
immediately and lives in this space forever; see also Wang \cite{W97}.
From Lemma~\ref{l3.2} we know that for any $r>0$ the family
$\{a_{i_j}: j=1,\cdots, m^\psi(r)\}$ under the regular
conditional probability $\P\{\,\cdot\,| m^\psi(r)\}$ are i.i.d.\
random variables with distribution $\<1,\mu\>^{-1} \mu(dx)$. This
gives an intuitive description of the locations $\{a_{i_j}:
j=1,\cdots, m^\psi (r)\}$ of the ``ancestors'' at the initial
time of $X_r$. By Lemma~\ref{l3.4}, each $\{w_{i_j} (a_{i_j},t):
t\ge r\}$ under $\P\{\,\cdot\,|{\cal G}_r\}$ is a Feller branching
diffusion with branching rate $\{\sigma(x(a_{i_j},t)): t\ge r\}$.
Then we have
 \beqlb\label{3.14}
dw_{i_j}(a_{i_j},t)
=
\sqrt{\sigma(x(a_{i_j},t)) w_{i_j}(a_{i_j},t)} dB_{i_j}(r,t),
\quad t\ge r,
 \eeqlb
for a Brownian motion $\{B_{i_j}(r,t): t\ge r\}$ stopped at
$\tau_0(w_{i_j}(a_{i_j}))$. However, under any enumeration the whole
excursion, $\{w_{i_j} (a_{i_j},t): t\ge 0\}$ is not a Feller
branching diffusion, otherwise the initial condition $w_{i_j}
(a_{i_j},0) =0$ would imply $w_{i_j}(a_{i_j},t)=0$ for all $t\ge 0$.
Therefore, the constructions (\ref{3.9}) and (\ref{3.12}) of the
SDSM are essentially different. Indeed, the purely atomic version of
the SDSM with a general initial state can only be constructed by
excursions, not usual Feller branching diffusions.


\section{SDSM with deterministic immigration}

\setcounter{equation}{0}

In this section, we construct some immigration processes by
one-dimensional excursions carried by stochastic flows. To simplify
the discussion, we assume the branching density is a constant
$\sigma >0$. Suppose that $m\in M(\IR)$ satisfies $\<1,m\>
>0$. Let $\L$ be given by (\ref{3.1}) and define
 \beqlb\label{4.1}
\J F(\nu)
=
\L F(\nu) + \int_{\IR} \frac{\delta F(\nu)}{\delta\nu(x)}m(dx),
\quad \nu\in M(\IR).
 \eeqlb
Setting ${\cal D}(\J) = {\cal D}(\L)$, we shall see that the
$(\J,{\cal D}(\J))$-martingale problem is equivalent with the one
given by (\ref{1.3}) and (\ref{1.4}) with deterministic
immigration rate $q(\cdot,\cdot) \equiv 1$.

Suppose that $(\itOmega,{\cal F},\P)$ is a complete standard
probability space on which we have: (i) a white noise $W(ds,dy)$
on $[0,\infty)\times \IR$ based on the Lebesgue measure; (ii) a
sequence of independent $\sigma$-branching diffusions $\{\xi_i(t):
t\ge0\}$ with $\xi_i(0) \ge0$ $(i=1,2,\cdots)$; (iii) a Poisson
random measure $N(ds,da,dw)$ with intensity $ds m(da) \Q_\kappa
(dw)$ on $[0,\infty) \times \IR\times \W_0$, where $\Q_\kappa$
denotes the excursion law of the $\sigma$-branching diffusion. We
assume that $\sum_{i=1} ^\infty \xi_i(0) < \infty$ and that
$\{W(ds,dy)\}$, $\{\xi_i(t)\}$ and $\{N(ds,da,dw)\}$ are
independent of each other. Given $(r,a)\in [0,\infty) \times \IR$,
let $\{x(r,a,t): t\ge r\}$ denote the unique solution of
 \beqlb\label{4.2}
x(t) = a + \int_r^t\int_{\IR} h(y-x(s))W(ds,dy), \quad t\ge r.
 \eeqlb
For $t\ge0$ let ${\cal G}_t$ be the $\sigma$-algebra generated by
all $\P$-null sets and the families of random variables
 \beqlb\label{4.3}
\{W([0,s]\times B), \xi_i(s): 0\le s\le t; B\in{\cal B}(\IR),
i=1,2,\cdots\}
 \eeqlb
and
 \beqlb\label{4.4}
\{N(J\times A): J\in{\cal B}([0,s]\times \IR);
A\in{\cal B}_{t-s}(\W_0); 0\le s\le t\}.
 \eeqlb
For a sequence $\{a_i\} \subset \IR$ let
 \beqlb\label{4.5}
Y_t = \sum_{i=1}^\infty \xi_i(t)\delta_{x(0,a_i,t)}
+ \int_0^t\int_{\IR}\int_{\W_0} w(t-s)
\delta_{x(s,a,t)}N(ds,da,dw),
\quad t\ge0.
 \eeqlb
(Here and in the sequel we make the convention that $\int_0^t =
\int_{(0,t]}$.)

We shall prove that $\{Y_t: t\ge0\}$ is a.s.\ continuous and
solves the $(\J,{\cal D}(\J))$-martingale problem relative to the
filtration $({\cal G}_t)_{t\ge0}$. Let $\W_{1/n} = \{w\in \W_0:
\tau_0(w)>1/n\}$ and recall from (\ref{2.8}) that $\Q_\kappa
(\W_{1/n}) = 2n/\sigma$. To prove the continuity of $\{Y_t: t\ge0\}$ we
consider the following approximating sequence:
 \beqlb\label{4.6}
Y^{(n)}_t = \sum_{i=1}^\infty \xi_i(t)\delta_{x(0,a_i,t)}
+ \int_0^t\int_{\IR}\int_{\W_{1/n}} w(t-s)\delta_{x(s,a,t)}
N(ds,da,dw),
\quad t\ge 0.
 \eeqlb

\blemma\label{l4.1} Both $\{Y_t: t\ge0\}$ and $\{Y^{(n)}_t: t\ge
0\}$ are a.s.\ continuous, and for any $T>0$ and $\phi\in
C(\IR)^+$ we have a.s.\ $\{\<\phi,Y^{(n)}_t\>: 0\le t\le T\}$
converges to $\{\<\phi,Y_t\>: 0\le t\le T\}$ increasingly and
uniformly as $n\to\infty$. \elemma

\noindent{\it Proof.} Let $N_1(ds,dw)$ denote the image of
$N(ds,da,dw)$ under the mapping $(s,a,w)\mapsto (s,w)$. Then
$N_1(ds,dw)$ is a Poisson random measure on $[0,\infty)\times
\W_0$ with intensity $\<m,1\> ds\Q_\kappa(dw)$ and is independent
of $\{\xi_i(t): t\ge0; i=1,2,\cdots\}$. Note that
 \beqnn
\<1,Y_t\>
=
\sum_{i=1}^\infty \xi_i(t)
+ \int_0^t\int_{\W_0} w(t-s) N_1(ds,dw),
\quad t\ge0.
 \eeqnn
By Pitman and Yor \cite[Theorem~4.1]{PY82}, $\{\<1,Y_t\>: t\ge0\}$ is a
diffusion process with generator $2^{-1}\sigma xd^2/dx^2 + \<1,m\>
d/dx$. Let $\itOmega_1 \in{\cal F}$ be a set with full
$\P$-measure such that $\{\<1,Y_t (\omega)\>: t\ge0\}$ is
continuous and $N(\omega, [0,n] \times \IR\times \W_{1/n})<\infty$ for
all $n\ge 1$ and $\omega\in \itOmega_1$. For $\omega\in
\itOmega_1$ and $\phi\in C(\IR)^+$, we have that $\{\<\phi,
Y^{(n)}_t (\omega)\>: t\ge0\}$ is continuous and converges to
$\{\<\phi,Y_t(\omega)\>: t\ge0\}$ increasingly as $n\to\infty$.
Then $\{\<\phi,Y_t(\omega)\>: t\ge0\}$ is lower semi-continuous.
The same reasoning shows that
 \beqnn
\<\|\phi\|-\phi,Y_t(\omega)\>
= \|\phi\|\<1,Y_t(\omega)\>
- \<\phi,Y_t(\omega)\>,
\quad t\ge 0,
 \eeqnn
is also lower semi-continuous. Since $\{\<1,Y_t(\omega)\>:
t\ge0\}$ is continuous, we conclude that $\{\<\phi,Y_t(\omega)\>:
t\ge0\}$ is continuous, giving the desired results. \qed

To show $\{Y_t: t\ge 0\}$ is a solution of the $(\J,{\cal
D}(\J))$-martingale problem relative to the filtration $({\cal
G}_t)_{t\ge0}$, we consider another approximating sequence
$\{X^{(n)}_t: t\ge 0\}$ defined by
 \beqlb\label{4.7}
X^{(n)}_t = \sum_{i=1}^\infty \xi_i(t)\delta_{x(0,a_i,t)}
+ \int_0^t\int_{\IR}\int_{\W_{1/n}} w(t-s+1/n)\delta_{x(s,a,t)}
N(ds,da,dw).
 \eeqlb
Note that the first parts of the excursions with length $1/n$
were chopped off in taking the second summation in (\ref{4.7}).

\blemma\label{l4.2} The process $\{X^{(n)}_t: t\ge0\}$ is a.s.\
c\`adl\`ag and for any $T>0$ and $\phi\in C(\IR)$ we have a.s.\
$\{\<\phi, X^{(n)}_t\>: 0\le t\le T\}$ converges to
$\{\<\phi,Y_t\>: 0\le t\le T\}$ uniformly as $n\to\infty$. \elemma

\noindent{\it Proof.} As in the proof of Lemma~\ref{l4.1}, we may
assume $\xi_i(\cdot) \equiv 0$ for all $i\ge 1$. The first
assertion holds since $N([0,n] \times \IR \times \W_{1/n}) <\infty$
a.s.\ for all $n\ge 1$. Clearly, we have a.s. $\<1, X^{(n)}_t\>
\le \<1, Y_{t+1/n}\>$ simultaneously for all $t\ge0$. By
Lemma~\ref{l4.1}, there is a set $\itOmega_2\in{\cal F}$ with full
$\P$-measure such that $\{\<1,Y^{(n)}_t\>: 0\le t\le T\}$
converges to $\{\<1,Y_t\>: 0\le t\le T\}$ uniformly as
$n\to\infty$ for all $T>0$ and $\omega\in \itOmega_2$. Fix $T>0$
and $\omega\in \itOmega_2$. For $\varepsilon >0$ let $m(\omega)\ge
1$ be an integer such that
 \beqnn
\<1,Y_t(\omega)\> - \<1,Y^{(n)}_t(\omega)\> < \epsilon,
\quad 0\le t\le T+1,
 \eeqnn
for $n\ge m(\omega)$ or, equivalently,
 \beqnn
\int_0^t\int_{\IR}\int_{\W_0\setminus \W_{1/m(\omega)}} w(t-s)
N(\omega,ds,da,dw)
< \epsilon,
\quad 0\le t\le T+1.
 \eeqnn
Since $N(\omega,[0,T]\times \IR \times \W_{1/m(\omega)})<\infty$,
there is an integer $M(\omega)\ge m(\omega)$ such that, for $n\ge
M(\omega)$,
 \beqnn
\int_0^t\int_{\IR}\int_{\W_{1/m(\omega)}} |w(t-s+1/n)-w(t-s)|
N(\omega,ds,da,dw)
< \epsilon,
\quad 0\le t\le T.
 \eeqnn
Then for $n\ge M(\omega)$ and $\phi\in C(\IR)$ we have
 \beqnn
& &|\<\phi,X^{(n)}_t(\omega)\> - \<\phi,Y_t(\omega)\>|   \\
&\le&
\int_0^t\int_{\IR}\int_{\W_{1/m(\omega)}}
\|\phi\||w(t-s+1/n)-w(t-s)| N(\omega,ds,da,dw)    \\
& &+ \int_0^t\int_{\IR}\int_{\W_0\setminus \W_{1/m(\omega)}}
\|\phi\||w(t-s+1/n)+w(t-s)| N(\omega,ds,da,dw)    \\
&<& 3\|\phi\|\epsilon
 \eeqnn
for $0\le t\le T$. That is, $\{\<\phi, X^{(n)}_t(\omega)\>: 0\le
t\le T\}$ converges to $\{\<\phi, Y_t(\omega)\>: 0\le t\le T\}$
uniformly as $n\to\infty$. \qed

We can easily pick out $\beta$-branching diffusions in the process
$\{X^{(n)}_t: t\ge0\}$. Let $N_n(ds,da,dw)$ denote the restriction
of $N(ds,da,dw)$ to $[0,\infty)\times \IR\times \W_{1/n}$. For $t\ge0$
let $\eta_n(t) = N([0,t]\times \IR \times \W_{1/n})$ and let ${\cal
G}_{n\!,t}$ be the $\sigma$-algebra generated by the families
(\ref{4.3}) and
 \beqlb\label{4.8}
\{N(J\times A): J\in{\cal B}([0,s]\times \IR);
A\in{\cal B}_{t-s+1/n}(\W_0); 0\le s\le t\}.
 \eeqlb
Clearly, we can a.s.\ arrange the atoms of $N_n(ds,da,dw)$ into
a sequence $\{(r_j,b_j,w_j): j=1,2,\cdots,\}$ so that
$0<r_1< r_2< \cdots$. With this ordering we have

\blemma\label{l4.3} For any $r\ge 0$, the family $\{\xi_i(t+r),
w_j(t+r-r_j+1/n): t\ge 0; i=1,2,\cdots; j=1,\cdots,\eta_n(r)\}$
under $\P\{\,\cdot\,|{\cal G}_{n\!,r}\}$ are independent
$\sigma$-branching diffusions relative to $({\cal
G}_{n\!,r+t})_{t\ge0}$. \elemma

\noindent{\it Proof.} For $r\ge0$ let ${\cal F}_{n\!,r}$ be the
$\sigma$-algebra generated by
 \beqnn
\{N_n(J\times A): J\in{\cal B}([0,r]\times \IR);
A\in{\cal B}_{1/n} (\W_0) \cap \W_{1/n}\}.
 \eeqnn
By Theorem~\ref{t2.1} and Lemma~\ref{l3.2}, $\{w_j(t+1/n): t\ge 0;
j=1,\cdots,\eta_n(r)\}$ under $\P\{\,\cdot\,|{\cal F}_{n\!,r}\}$
are independent $\sigma$-branching diffusions. Clearly, the same
assertion is true for $\{w_j(t+r-r_j+1/n): t\ge 0; j=1, \cdots,
\eta_n(r)\}$ under $\P\{\,\cdot\,|{\cal F}_{n\!,r}^\prime\}$,
where ${\cal F}_{n\!,r}^\prime$ is the $\sigma$-algebra generated
by
 \beqnn
\{N_n(J\times A): J\in{\cal B}([0,s]\times \IR);
A\in{\cal B}_{r-s+1/n} (\W_0) \cap \W_{1/n}; 0\le s\le r\}.
 \eeqnn
Then the result follows from the independence of $\{\xi_n(t)\}$,
$\{W(ds,dy)\}$ and $\{N(ds,da,dw)\}$. \qed

\blemma\label{l4.4} The process $\{X^{(n)}_t: t\ge0\}$ relative to
$({\cal G}_{n\!,t})_{t\ge0}$ is a strong Markov process generated
by the closure of $(\J_n,{\cal D}(\J_n))$, where
 \beqlb\label{4.9}
\J_n F(\nu)
=
\L F(\nu) + \int_{\IR} m(dx)\int_0^\infty [F(\nu+y\delta_x)
- F(\nu)]\kappa_{1/n}(dy),
\quad \nu\in M(\IR),
 \eeqlb
for $F\in {\cal D}(\J_n) = {\cal D}(\L)$. \elemma

\noindent{\it Proof.} Clearly, each $r_k$ is a stopping times and
for $0\le t<r_{k+1}-r_k$ we have
 \beqnn
X^{(n)}_{t+r_k} =
\sum_{i=1}^\infty \xi_i(t+r_k)\delta_{x_i(0,a_i,r_k+t)}
+ \sum_{j=1}^k w_j(t+r_k-r_j+1/n) \delta_{x(r_j,b_j,r_k+t)}.
 \eeqnn
Since $r\ge0$ in Lemma~\ref{l4.3} was arbitrary, $\{\xi_i(t+r_k),
w_j(t+r_k-r_j+1/n): t\ge 0; i = 1,2, \cdots; j= 1, \cdots, k\}$
under $\P\{\,\cdot\,| {\cal G}_{n\!,r_k}\}$ are independent
$\sigma$-branching diffusions relative to $({\cal G}_{n\!,r_k+t})
_{t\ge0}$. By the independence of $\{\xi_n(t)\}$, $\{W(ds,dy)\}$
and $\{N(ds,da,dw)\}$ we may apply Theorem~\ref{t3.3} to get that
$\{X^{(n)}_{t+r_k}: 0\le t<r_{k+1}-r_k\}$ under $\P\{\,\cdot\,|
{\cal G}_{n\!,r_k}\}$ is a (killed) diffusion process relative to
$({\cal G}_{n\!,r_k+t})_{t\ge0}$ with generator $\L -
2n\<1,m\>/\sigma$. Observe also that
 \beqnn
\P\{F(X^{(n)}_{r_k})|{\cal G}_{n\!,r_k-}\}
=
\frac{\sigma}{2n\<1,m\>}\int_{\IR}m(dx)\int_0^\infty
F(X^{(n)}_{r_k-}+y\delta_x)\kappa_{1/n}(dy)
 \eeqnn
for any $F\in C(M(\IR))$. Then $\{X^{(n)}_t: t\ge 0\}$ relative to
$({\cal G} _{n\!,t})_{t\ge0}$ is a strong Markov process generated
by the closure of $(\J_n,{\cal D}(\J_n))$. \qed

\btheorem\label{t4.1} The process $\{Y_t: t\ge0\}$ constructed by
(\ref{4.5}) is a.s.\ continuous and solves the $(\J,{\cal
D}(\J))$-martingale problem relative to the filtration $({\cal
G}_t)_{t\ge0}$. \etheorem

\noindent{\it Proof.} As observed in the proof of
Lemma~\ref{l4.1}, $\{\<1,Y_t(\omega)\>: t\ge0\}$ is a diffusion
process with generator $2^{-1}\sigma xd^2/dx^2 + \<1,m\> d/dx$ and
initial value $\<1,\mu\>$. Then $\P\{\<1,Y_t\>^n\}$ is a locally
bounded function of $t\ge0$ for every $n\ge1$. For
$F=F_{f,\{\phi_i\}}$ given by (\ref{3.2}), Lemma~\ref{l4.4}
asserts that
 \beqnn
F(X^{(n)}_t) - F(X^{(n)}_0) - \int_0^t \J_n F(X^{(n)}_s)ds,
\quad t\ge0,
 \eeqnn
is a martingale relative to $({\cal G}_{n\!,t})_{t\ge0}$. Since
$({\cal G}_t)_{t\ge0}$ is smaller than $({\cal G}_{n\!,t})_{t\ge0}$,
for any $t\ge r\ge0$ and any $G\in\mbox{bp}{\cal G}_r$ we have
 \beqlb\label{4.10}
\E\bigg\{G\bigg[F(X^{(n)}_t) - F(X^{(n)}_r)
- \int_r^t \J_n F(X^{(n)}_s)ds\bigg]\bigg\}
= 0.
 \eeqlb
By (\ref{2.5}) and (\ref{4.9}) it is not hard to check that $\J_n
F(\nu) \to \J F(\nu)$ uniformly on the set $\{\nu\in M(\IR):
\<1,\nu\>\le a\}$ as $n\to \infty$ for each $a\ge0$. By
Lemma~\ref{l4.2}, letting $n\to\infty$ in (\ref{4.10}) we get
 \beqlb\label{4.11}
\E\bigg\{G\bigg[F(Y_t) - F(Y_r)
- \int_r^t \J F(Y_s)ds\bigg]\bigg\}
= 0.
 \eeqlb
That is,
 \beqnn
F(Y_t) - F(Y_0) - \int_0^t \J F(Y_s)ds,
\quad t\ge0,
 \eeqnn
is a martingale relative to $({\cal G}_t)_{t\ge0}$. Then the
desired result follows by an approximation of an arbitrary
$F\in {\cal D}(\J)$. \qed

We can also pick out $\sigma$-branching diffusions in the process
$\{Y_t: t\ge0\}$ defined by (\ref{4.5}). For $r>0$ we can a.s.\
enumerate the atoms $(s,a,w)$ of $N(ds,da,dw)$ satisfying $0< s
<r$ and $w(r-s)>0$ into a sequence $\{(r_j,b_j,w_j): j=1,2,
\cdots\}$ so that $r_j <r_{j+1}$ for all $j\ge1$. This enumeration
gives the following

\btheorem\label{t4.2} For any $r>0$, the sequence $\{\xi_i(t+r),
w_j(t+r-r_j): t\ge 0; i=1,2,\cdots; j=1,2,\cdots\}$ under
$\P\{\,\cdot\,|{\cal G}_r\}$ are independent $\sigma$-branching
diffusions relative to $({\cal G}_{r+t})_{t\ge0}$. \etheorem

\noindent{\it Proof.} Clearly, for any integer $n>1/r$, we have
$r_j<r-1/n$ if and only if $j\le \eta_n (r-1/n)$. As in the proof
of Lemma~\ref{l4.3} we see that, the sequence $\{w_j(t+r-r_j):
t\ge 0; j=1, \cdots, \eta_n(r-1/n)\}$ under $\P\{\,\cdot\,|{\cal
F}_{n\!,r} ^{\prime\prime}\}$ are independent $\sigma$-branching
diffusions, where ${\cal F}_{n\!,r}^{\prime\prime}$ is the
$\sigma$-algebra generated by
 \beqnn
\{N_n(J\times A): J\in{\cal B}([0,s]\times \IR);
A\in{\cal B}_{r-s} (\W_0) \cap \W_{1/n}; 0\le s\le r\}.
 \eeqnn
By the independence of $\{\xi_i(t)\}$, $\{W(ds,dy)\}$ and
$\{N(ds,da,dw)\}$ we have that $\{\xi_i(t+r), w_j(t+r-r_j): t\ge
0; i=1,2,\cdots; j=1, \cdots, \eta_n(r-1/n)\}$ under
$\P\{\,\cdot\,| {\cal G}_r\}$ are independent $\sigma$-branching
diffusions. Since $\eta_n(r-1/n) \to \infty$ as $n\to \infty$,
we have the desired result. \qed

Now let us consider an arbitrary initial state $\mu\in M(\IR)$.
Suppose on the complete standard probability space $(\itOmega,
{\cal F}, \P)$ have: (i) a white noise $W(ds,dy)$ on $[0,\infty)
\times \IR$ based on the Lebesgue measure; (ii) a Poisson random
measure $N_0(da,dw)$ on $\IR\times \W_0$ with intensity $\mu(da)
\Q_\kappa(dw)$; and (iii) a Poisson random measure $N(ds,da,dw)$
on $[0,\infty) \times \IR\times \W_0$ with intensity $ds
m(da)\Q_\kappa(dw)$. We assume that $\{W(ds,dy)\}$,
$\{N_0(da,dw)\}$ and $\{N(ds,da,dw)\}$ are independent of each
other. For $t\ge0$ let ${\cal G}_t$ be the $\sigma$-algebra
generated by all $\P$-null sets and the families of random
variables
 \beqlb\label{4.12}
\{W([0,s]\times B), N_0(F\times A): F\in{\cal B}(\IR);
A\in{\cal B}_t(\W_0); B\in{\cal B}(\IR); 0\le s\le t\}
 \eeqlb
and
 \beqlb\label{4.13}
\{N(I\times B\times A): I\in{\cal B}([0,s]); B\in{\cal B}(\IR);
A\in{\cal B}_{t-s}(\W_0); 0\le s\le t\}.
 \eeqlb
Given $(r,a)\in [0,\infty) \times \IR$, let $\{x(r,a,t): t\ge r\}$
denote the unique solution of (\ref{4.2}). Let $Y_0=\mu$ and for
$t>0$ let
 \beqlb\label{4.14}
Y_t = \int_{\IR}\int_{\W_0} w(t) \delta_{x(0,a,t)}N_0(da,dw)
+ \int_0^t\int_{\IR}\int_{\W_0} w(t-s) \delta_{x(s,a,t)} N(ds,da,dw).
 \eeqlb

\btheorem\label{t4.3} The process $\{Y_t: t\ge0\}$ defined above
is a.s.\ continuous and solves the $(\J,{\cal D}(\J))$-martingale
problem relative to the filtration $({\cal G}_t)_{t\ge0}$.
\etheorem

\noindent{\it Proof.} By Theorem~\ref{t3.4}, the first term on the
right hand side of (\ref{4.14}) is a.s.\ continuous and converges
to $\mu$ as $t\to0$. By Lemma~\ref{l4.1} the second term is also
a.s.\ continuous. Thus $\{Y_t: t\ge0\}$ is a.s.\ continuous. For
$r>0$ let $\{(a_i,u_i): i=1,\cdots, m(r)\}$ be set of atoms of
$N_0(da,dw)$ satisfying $u_i(r)>0$ and be arranged so that $u_1(r)
< \cdots < u_{m(r)}(r)$.  Let $\{(r_j,b_j,w_j): j=1,2, \cdots\}$
be the set of atoms of $N(ds,da,dw)$ satisfying $0< r_j < r$ and
$w_j(r-r_j)>0$ and be arranged so that $r_j <r_{j+1}$ for all
$j\ge1$. By Lemma~\ref{l3.4}, Theorem~\ref{t4.2} and the
independence assumption, $\{u_i(t+r), w_j(t+r-r_i): t\ge 0;
i=1,\cdots,m(r); j=1,2, \cdots\}$ under $\P\{\,\cdot\,|{\cal
G}_r\}$ are independent $\sigma$-branching diffusions relative to
$({\cal G}_{r+t})_{t\ge0}$. By Theorem~\ref{t4.1} and the property
of independent increments of $W(ds,dy)$ and $N(ds,da,dw)$, the
continuous process $\{Y_{t+r}: t\ge0\}$ under $\P\{\,\cdot\,|{\cal
G}_r\}$ is a solution of the $(\J,{\cal D}(\J))$-martingale
problem relative to the filtration $({\cal G}_{t+r})_{t\ge0}$.
Since $r>0$ was arbitrary in the above reasoning, we have the
desired result. \qed

\btheorem\label{t4.4} The $(\J,{\cal D}(\J))$-martingale problem has
a unique solution. \etheorem

\noindent{\it Proof.} By Theorem~\ref{t4.3}, the $(\J,{\cal
D}(\J))$-martingale problem has a solution. The approach to the
uniqueness is similar to that in Dawson {\it et al} \cite{DLW01},
so we only provide an outline. For $f\in C^2(\IR^n)$ let
 \beqnn
G^nf(x)
=
\frac{1}{2}\rho(0)\sum_{i=1}^n \frac{\partial^2}{\partial x_i^2}f(x)
+ \frac{1}{2} \sum_{i,j=1, i\neq j}^n
\rho(x_i-x_j)\frac{\partial^2}{\partial x_i \partial x_j}f(x),
\quad x\in \IR^n.
 \eeqnn
Define $\itPhi_{ij}f \in C(\IR^{n-1})$ by
 \beqnn
\itPhi_{ij}f(x_1,\cdots,x_{n-1})
=
\sigma(x_{n-1})f(x_1,\cdots,x_{n-1},\cdots,x_{n-1},\cdots,x_{n-2}),
 \eeqnn
where $x_{n-1}\in \IR$ is in the places of the $i$th and the $j$th
variables of $f$ on the right hand side, and define $\itPsi_i f\in
C^2 (\IR^{n-1})$ by
 \beqnn
\itPsi_if(x_1,\cdots,x_{n-1})
=
\int_{\IR} f(x_1,\cdots,x_{i-1},x,x_i,\cdots,x_{n-1}) m(dx),
 \eeqnn
where $x\in \IR$ is the $i$th variable of $f$ on the right hand
side. It is not hard to show that
 \beqnn
\J F_{n,f}(\nu)
=
F_{n,G^nf}(\nu)
+ \frac{1}{2} \sum_{i,j=1, i\neq j}^n F_{n-1,\itPhi_{ij}f}(\nu)
+ \sum_{i=1}^n F_{n-1,\itPsi_if}(\nu),
\quad \nu\in M(\IR).
 \eeqnn
Write $F_\nu(n,f) = F_{n,f}(\nu)$ and let
 \beqlb\label{4.15}
\J^* F_\nu(n,f)
&=&
F_\nu(n,G^nf)
+ \,\frac{1}{2} \sum_{i,j=1, i\neq j}^n
[F_\nu(n-1,\itPhi_{ij}f) - F_\nu(n,f)]     \nonumber  \\
& &
+ \,\sum_{i=1}^n [F_\nu(n-1,\itPsi_if) - F_\nu(n,f)].
 \eeqlb
Then we have
 \beqlb\label{4.16}
\J F_{n,f}(\nu)
=
\J^* F_\nu(n,f) + \frac{1}{2} n(n+1) F_\nu(n,f).
 \eeqlb
Guided by (\ref{4.15}) we can construct a Markov process $\{(M_t,
F_t): t\ge 0\}$ with initial value $(M_0,F_0) = (n,f)$ and
generator $\J^*$. Based on (\ref{4.16}) one can prove that if
$\{Y_t: t\ge 0\}$ is a solution of the $(\J,{\cal
D}(\J))$-martingale problem with $Y_0=\mu$, then
 \beqnn
\E_\mu\{\<f,Y_t^n\>\}
=
\E_{(n,f)} \bigg[\<F_t, \mu^{M_t}\>
\exp\bigg\{\frac{1}{2}\int_0^t M_s(M_s+1)ds \bigg\}\bigg],
\quad t\ge0;
 \eeqnn
see \cite[p.195]{EK86}. This duality determines the
one-dimensional distributions of $\{Y_t: t\ge0\}$ uniquely, and
hence the conclusion follows by \cite[p.184]{EK86}. \qed

By the uniqueness of solution of the $(\J$, ${\cal
D}(\J))$-martingale problem, the immigration process constructed
by (\ref{4.14}) is a diffusion. From this construction we know
that the immigration SDSM started with any initial state actually
lives in the space purely atomic measures. The next theorem, which
can be proved similarly as Theorem~\ref{t3.2}, gives a useful
alternate characterization of the immigration SDSM.

\btheorem\label{t4.5} A continuous $M(\IR)$-valued process $\{Y_t:
t\ge0\}$ is a solution of the $(\J$, ${\cal D}(\J))$-martingale
problem if and only if for each $\phi\in C^2(\IR)$,
 \beqlb\label{4.17}
M_t(\phi)
:=
\<\phi, Y_t\> - \<\phi, Y_0\> - \<\phi, m\> t
- \frac{1}{2} \rho(0)\int_0^t \<\phi^{\prime\prime}, Y_s\> ds,
\quad t\ge0,
 \eeqlb
is a martingale with quadratic variation process
 \beqlb\label{4.18}
\<M(\phi)\>_t
=
\int_0^t\<\sigma\phi^2, Y_s\>ds
+ \int_0^t ds\int_{\IR} \<h(z - \cdot) \phi^\prime, Y_s\>^2 dz.
 \eeqlb
\etheorem

Under the condition of Theorem~\ref{t4.5}, the martingales
$\{M_t(\phi): t\ge 0\}$ defined by (\ref{4.17}) and (\ref{4.18})
form a system which is linear in $\phi \in C^2 (\IR)$. Following
the method of Walsh \cite{W86}, we can define the stochastic integral
 \beqnn \int_0^t\int_{\IR}\phi(s,x)M(ds,dx), \qquad t\ge 0,  \eeqnn
if both $\phi(s,x)$ and $\phi^\prime_x(s,x)$ are continuous on
$[0,\infty) \times \IR$. By a standard argument we get the
following

\btheorem\label{t4.6} In the situation described above, for any
$t\ge0$ and $\phi \in C^1 (\IR)$ we have a.s.
 \beqnn
\<\phi, Y_t\>
=
\<P_t\phi, Y_0\> + \int_0^t \<P_{t-s}\phi,m\>ds
+ \int_0^t\int_{\IR} P_{t-s}\phi(x) M(ds,dx),
 \eeqnn
where $(P_t)_{t\ge0}$ is the semigroup of the Brownian motion
generated by $2^{-1} \rho(0) d^2/dx^2$. \etheorem


\section{SDSM with interactive immigration}

\setcounter{equation}{0}

In this section, we construct a diffusion solution of the
martingale problem given by (\ref{1.3}) and (\ref{1.4}) with a
general interactive immigration rate. This is done by solving a
stochastic equation carried by a stochastic flow and driven by
Poisson processes of excursions.

Let $\sigma >0$ be a constant and let $m$ be a non-trivial
$\sigma$-finite Borel measure on $\IR$. Suppose we have on a
complete standard probability space $(\itOmega,{\cal F},\P)$ the
following: (i) a white noise $W(ds,dy)$ on $[0,\infty)\times \IR$
based on the Lebesgue measure; (ii) a sequence of independent
$\sigma$-branching diffusions $\{\xi_i(t): t\ge0\}$ with $\xi_i(0)
\ge0$ $(i= 1,2,\cdots)$; and (iii) a Poisson random measure
$N(ds,da,du,dw)$ on $[0,\infty) \times \IR\times [0,\infty) \times
\W_0$ with intensity $ds m(da)du \Q_\kappa (dw)$, where
$\Q_\kappa$ denotes the excursion law of the $\sigma$-branching
diffusion. We assume that $\sum_{i=1}^\infty \xi_i(0) <\infty$ and
that $\{W(ds,dy)\}$, $\{\xi_i(t)\}$ and $\{N(ds,da,du,dw)\}$ are
independent of each other. For $t\ge0$ let ${\cal G}_t$ be the
$\sigma$-algebra generated by all $\P$-null sets and the families
of random variables (\ref{4.3}) and
 \beqlb\label{5.1}
\{N(J\times A): J\in{\cal B}([0,s]\times \IR\times [0,\infty));
A\in{\cal B}_{t-s}(\W_0); 0\le s\le t\}.
 \eeqlb
Let ${\cal P}$ be the $\sigma$-algebra on $[0,\infty)\times
\IR\times \itOmega$ generated by functions of the form
 \beqlb\label{5.2}
g(s,x,\omega) = \eta_0(x,\omega)1_{\{0\}}(s)
+ \sum_{i=0}^\infty \eta_i(x,\omega)1_{(r_i,r_{i+1}]}(s),
 \eeqlb
where $0=r_0<r_1<r_2<\ldots$ and $\eta_i(\cdot,\cdot)$ is ${\cal
B}(\IR)\times {\cal G}_{r_i}$-measurable. We say a function on
$[0,\infty)\times \IR\times \itOmega$ is {\it predictable} if it is
${\cal P}$-measurable.

We first construct an immigration process with purely atomic
initial state and predictable immigration rate. Suppose that
$q(\cdot,\cdot,\cdot)$ is a non-negative predictable function on
$[0,\infty)\times \IR\times \itOmega$ such that $\E\{\<q(t,\cdot),
m\>^2\}$ is locally bounded in $t\ge0$. For a sequence
$\{a_i\}\subset \IR$ let
 \beqlb\label{5.3}
Y_t = \sum_{i=1}^\infty \xi_i(t)\delta_{x(0,a_i,t)}
+ \int_0^t\int_{\IR}\int_0^{q(s,a)}\int_{\W_0} w(t-s)
\delta_{x(s,a,t)}N(ds,da,du,dw),
\quad t\ge0.
 \eeqlb

\btheorem\label{t5.1} The process $\{Y_t: t\ge0\}$ defined by
(\ref{5.3}) has a continuous modification. For this
modification and each $\phi\in C^2(\IR)$,
 \beqlb\label{5.4}
M_t(\phi)
:=
\<\phi, Y_t\> - \<\phi, Y_0\>
- \frac{1}{2} \rho(0)\int_0^t \<\phi^{\prime\prime}, Y_s\> ds
- \int_0^t \<q(s,\cdot)\phi,m\>ds,
\quad t\ge0,
 \eeqlb
is a continuous martingale relative to the filtration $({\cal
G}_t)_{t\ge0}$ with quadratic variation process
 \beqlb\label{5.5}
\<M(\phi)\>_t
=
\int_0^t\<\sigma\phi^2, Y_s\>ds
+ \int_0^t ds\int_{\IR} \<h(z - \cdot) \phi^\prime, Y_s\>^2 dz.
 \eeqlb
\etheorem

\noindent{\it Proof.} {\it Step~1)} Suppose that $q(s,x,\omega)
\equiv q(x)$ for a function $q\in L^1(\IR,m)$. Let
$N_q(ds,da,du,dw)$ denote the restriction of $N(ds,da,du,dw)$ to
the set $\{(s,a,u,w): s\ge0; a\in\IR; 0\le u\le q(a); w\in \W_0\}$
and let $N_q(ds,da,dw)$ be the image of $N_q(ds,da,du,dw)$ under
the mapping $(s,a,u,w) \mapsto (s,a,w)$. Clearly, $N_q(ds,da,dw)$
is a Poisson measure on $[0,\infty)\times \IR\times \W_0$ with
intensity $dsq(a)m(da)\Q_\kappa(dw)$ and (\ref{5.3}) can be
rewritten as
 \beqnn
Y_t = \sum_{i=1}^\infty \xi_i(t)\delta_{x(0,a_i,t)}
+ \int_0^t\int_{\IR}\int_{\W_0} w(t-s)
\delta_{x(s,a,t)}N_q(ds,da,dw),
\quad t>0.
 \eeqnn
Then the results are reduced to those of Theorems~\ref{t4.1} and
\ref{t4.5}. {\it Step~2)} Suppose that $q(\omega,s,x)$ is of the
form (\ref{5.2}). Note that $\eta_i(x)$ is actually deterministic
under the conditional probability $\P\{\,\cdot\,|{\cal
G}_{r_i}\}$. By the last step and Theorems~\ref{t4.1} and
\ref{t4.5}, the results hold on each interval $[r_i,r_{i+1}]$ and
hence on $[0,\infty)$. {\it Step~3)} The case of a general
non-negative predictable function $q(\cdot,\cdot,\cdot)$ can be
proved by approximating arguments similar to those in Fu and Li
\cite{FL03} and Shiga \cite{S90}. \qed

Let us consider a stochastic equation with purely atomic initial
state. Suppose that $q(\cdot,\cdot)$ is a Borel function on
$M(\IR)\times \IR$ such that there is a constant $K$ such that
 \beqlb\label{5.6}
\<q(\nu,\cdot),m\>
\le
K(1 + \|\nu\|),
\quad \nu\in M(\IR),
 \eeqlb
and for each $R>0$ there is a constant $K_R>0$ such that
 \beqlb\label{5.7}
\<|q(\nu,\cdot) - q(\gamma,\cdot)|,m\>
\le
K_R \|\nu-\gamma\|
 \eeqlb
for $\nu$ and $\gamma\in M(\IR)$ satisfying $\<1,\nu\>\le R$ and
$\<1,\gamma\>\le R$, where $\|\cdot\|$ denotes the total variation.
For any sequence $\{a_i\}\subset\IR$, consider the stochastic
equation:
 \beqlb\label{5.8}
Y_t
= \sum_{i=1}^\infty \xi_i(t)\delta_{x(0,a_i,t)}
+ \int_0^t\int_{\IR}\int_0^{q(Y_s,a)}\int_{\W_0} w(t-s)
\delta_{x(s,a,t)}N(ds,da,du,dw),
\quad t\ge0.
 \eeqlb

\btheorem\label{t5.2} Under the above conditions, there is a
unique continuous solution $\{Y_t: t\ge0\}$ of (\ref{5.8}), which is
a diffusion process. Moreover, for each $\phi\in C^2(\IR)$,
 \beqlb\label{5.9}
M_t(\phi) = \<\phi,Y_t\> - \<\phi,Y_0\>
- \frac{1}{2}\int_0^t \<\phi^{\prime\prime},Y_s\> ds
- \int_0^t \<q(Y_s,\cdot)\phi,m\>ds,
\quad t\ge0,
 \eeqlb
is a continuous martingale relative to the filtration $({\cal
G}_t)_{t\ge0}$ with quadratic variation process
 \beqlb\label{5.10}
\<M(\phi)\>_t = \int_0^t\<\sigma\phi^2, Y_s\>ds +
\int_0^t ds\int_{\IR} \<h(z - \cdot) \phi^\prime, Y_s\>^2 dz,
\quad t\ge 0.
 \eeqlb
\etheorem

\noindent{\it Proof.} Based on Theorem~\ref{t5.1} and the results
in the last section, it can be proved by iteration arguments
similar to those in \cite{FL03} and \cite{S90} that
(\ref{5.8}) has a unique solution and (\ref{5.9}) is a continuous
martingale with quadratic variation process (\ref{5.10}). Let
$\mu= \sum_{i=1}^\infty \xi_i(0) \delta_{a_i}$ and let
$Q_t^q(\mu,\cdot)$ denote the distribution of $Y_t$ defined by
(\ref{5.8}). For any bounded $({\cal G}_t)$-stopping time
$\tau\ge0$, we can use the information from ${\cal G}_\tau$ to
enumerate the atoms $(s,a,u,w)$ of $N(ds,da,du,dw)$ satisfying $0<
s\le \tau$ and $w(\tau-s)>0$ into a sequence $\{(r_j,b_j,u_j,w_j):
j=1,2, \cdots\}$ so that $r_j \le r_{j+1}$ for all $j\ge1$. By the
strong Markov property of $\sigma$-branching diffusions and a
slight modification of the proof of Theorem~\ref{t4.2}, we can
show that $\{\xi_i(t+\tau), w_j(t+\tau-r_j): t\ge 0; i=1,2,\cdots;
j=1,2,\cdots\}$ under $\P\{\,\cdot\,|{\cal G}_\tau\}$ are
independent $\sigma$-branching diffusions relative to $({\cal
G}_{\tau+t})_{t\ge0}$. By the property of independent increments,
$W_\tau(ds,dy) := W(ds+\tau,dy)$ under $\P\{\,\cdot\,|{\cal
G}_\tau\}$ is a white noise on $[0,\infty) \times \IR$ based on
the Lebesgue measure. Similarly, $N_\tau(ds,da,du,dw) :=
N(ds+\tau,da,du,dw)$ under $\P\{\,\cdot\,|{\cal G}_\tau\}$ is a
Poisson random measure with intensity $dsm(da)du\Q_\kappa(dw)$.
Moreover, the families $\{\xi_i(t+\tau), w_j(t+\tau-r_j)\}$,
$\{W_\tau(ds,dy)\}$ and $\{N_\tau (ds,da,du,dw)\}$ under
$\P\{\,\cdot\,|{\cal G}_\tau\}$ are independent of each other.
Observe that
 \beqnn
Y_\tau
&=& \sum_{i=1}^\infty \xi_i(\tau)\delta_{x(0,a_i,\tau)}
+ \sum_{j=1}^\infty w_j(\tau-r_j)\delta_{x(r_j,b_j,\tau)}
 \eeqnn
and
 \beqnn
Y_{t+\tau}
&=& \sum_{i=1}^\infty \xi_i(t+\tau)\delta_{x(0,a_i,t+\tau)}
+ \sum_{j=1}^\infty w_j(t+\tau-r_j)\delta_{x(r_j,b_j,t+\tau)} \\
& & + \int_0^t\int_{\IR}\int_0^{q(Y_{s+\tau},a)}\int_{\W_0} w(t-s)
\delta_{x(s+\tau,a,t+\tau)}N_\tau(ds,da,du,dw).
 \eeqnn
By the uniqueness of solution of (\ref{5.8}), $Y_{t+\tau}$ under
$\P\{\,\cdot\,|{\cal G}_\tau\}$ has distribution $Q_t^q (Y_\tau,
\cdot)$, giving the strong Markov property of $\{Y_t: t\ge0\}$.
\qed

We now consider a stochastic equation with a general initial state
$\mu \in M(\IR)$. Suppose on the complete standard probability
space $(\itOmega, {\cal F}, \P)$ we have the following: (i) a
white noise $W(ds,dy)$ on $[0,\infty)\times \IR$ based on the
Lebesgue measure; (ii) a Poisson random measure $N_0(da,dw)$ on
$\IR\times \W_0$ with intensity $\mu(dx) \Q_\kappa(dw)$; and (iii)
a Poisson random measure $N(ds,da,du,dw)$ on $[0,\infty) \times
\IR \times [0,\infty) \times \W_0$ with intensity $ds m(da) du
\Q_\kappa(dw)$. We assume that $\{W(ds,dy)\}$, $\{N_0(da,dw)\}$
and $\{N(ds,da,dw)\}$ are independent of each other. For $t\ge0$
let ${\cal G}_t$ be the $\sigma$-algebra generated by all
$\P$-null sets and the families of random variables (\ref{4.12})
and (\ref{5.1}).

\btheorem\label{t5.3} Suppose that $q(\cdot,\cdot)$ is a Borel
function on $M(\IR)\times \IR$ satisfying (\ref{5.6}) and
(\ref{5.7}). Then the stochastic equation:
 \beqlb\label{5.11}
Y_t
&=& \int_{\IR}\int_{\W_0} w(t) \delta_{x(0,a,t)}
N_0(da,dw)  \nonumber \\
& & + \int_0^t\int_{\IR}\int_0^{q(Y_s,a)}\int_{\W_0} w(t-s)
\delta_{x(s,a,t)}N(ds,da,du,dw),
\quad t> 0,
 \eeqlb
has a unique continuous solution $\{Y_t: t>0\}$. If we set
$Y_0=\mu$, then $\{Y_t: t\ge0\}$ is a diffusion process and the
martingale characterization of Theorem~\ref{t5.2} holds. \etheorem

\noindent{\it Proof.} If $q(\cdot,\cdot,\cdot)$ is a non-negative
predictable function on $[0,\infty)\times \IR\times \itOmega$ such
that $\E\{m(q(t,\cdot))^2\}$ is locally bounded in $t\ge0$, it
can be proved in three steps as in the proof of Theorem~\ref{t5.1}
that the process $\{Y_t: t\ge0\}$ defined by $Y_0=\mu$ and
 \beqnn
Y_t
&=& \int_{\IR}\int_{\W_0} w(t) \delta_{x(0,a,t)}
N_0(da,dw)  \nonumber \\
& & + \int_0^t\int_{\IR}\int_0^{q(s,a)}\int_{\W_0} w(t-s)
\delta_{x(s,a,t)}N(ds,da,du,dw),
\quad t> 0,
 \eeqnn
has a continuous modification and the results of
Theorem~\ref{t5.1} also hold for this process. Then one can show
by iteration arguments that (\ref{5.11}) has a unique solution and
the martingale characterization of Theorem~\ref{t5.2} holds. The
strong Markov property of $\{Y_t: t\ge0\}$ can be derived as in
the proof of Theorem~\ref{t5.2}. \qed

The solution of (\ref{5.11}) can be regarded as immigration
processes associated with the SDSM with interactive immigration.
It is not hard to show that the generator of the diffusion process
$\{Y_t: t\ge0\}$ is given by
 \beqlb\label{5.12}
\J F(\nu) =
\L F(\nu) + \int_{\IR} q(\nu,x)\frac{\delta F(\nu)}
{\delta\nu(x)} m(dx),\quad \nu\in M(\IR),
 \eeqlb
where $\L$ is defined by (\ref{3.1}) and $q(\cdot,\cdot)$ is the
interactive immigration rate. Note that the Markov property of
$\{Y_t: t\ge0\}$ was obtained from the uniqueness of solution of
(\ref{5.8}). This application of the stochastic equation is
essential since the uniqueness of solution of the martingale
problem given by (\ref{5.9}) and (\ref{5.10}) still remains open;
see also Fu and Li \cite{FL03} and Shiga \cite{S90}.

\medskip

{\bf Acknowledgement.} We thank Hao Wang and a referee for their
helpful comments and suggestions on an earlier version of the
paper.



\noindent
\begin{thebibliography}{99}


\bibitem{D93}
Dawson, D.A.: {\it Measure-Valued Markov Processes}. In: Lect.
Notes. Math. {\bf  1541}, 1-260, Springer-Verlag, Berlin (1993).

\bibitem{DLW01}
Dawson, D.A.; Li, Z.H. and Wang, H.: {\it Superprocesses with
dependent spatial motion and general branching densities}. Elect.
J. Probab. {\bf 6} (2001), Paper No. 25, 1-33.

\bibitem{DVW00}
Dawson, D.A.; Vaillancourt, J. and Wang, H.: {\it Stochastic
partial differential equations for a class of measure-valued
branching diffusions in a random medium}. Ann. Inst. Henri
Poincare, Probabilit\'es and Statistiques {\bf 36} (2000),
167-180.

\bibitem{EK86}
Ethier, S.N. and Kurtz, T.G.: {\it Markov Processes:
Characterization and Convergence}. Wiley, New York (1986).

\bibitem{FL03}
Fu, Z.F. and Li, Z.H.: {\it Measure-valued diffusions and
stochastic equations with Poisson process}. Osaka J. Math.,
submitted (2003), ps and pdf files: {\tt
math.bnu.edu.cn/\~{}lizh}.

\bibitem{IW89}
Ikeda, N. and Watanabe, S.: {\it Stochastic Differential Equations
and Diffusion Processes}. North-Holland/Kodansha, Amsterdam/Tokyo
(1989).

\bibitem{KS88}
Konno, N. and Shiga, T.: {\it Stochastic partial differential
equations for some measure-valued diffusions}. Probab. Theory
Related Fields {\bf 79} (1988), 201-225.

\bibitem{L02}
Li, Z.H.: {\it Skew convolution semigroups and related immigration
processes}. Theory Probab. Appl. {\bf 46} (2002), 274-296.

\bibitem{LS95}
Li, Z.H. and Shiga T.: {\it Measure-valued branching diffusions:
immigrations, excursions and limit theorems}. J. Math. Kyoto Univ.
{\bf 35} (1995), 233-274.

\bibitem{PY82}
Pitman, J. and Yor, M.: {\it A decomposition of Bessel bridges}.
Z. Wahrsch. verw. Geb. {\bf 59} (1982), 425-457.

\bibitem{R89}
Reimers, M.: {\it One dimensional stochastic differential
equations and the branching measure diffusion}. Probab. Theory
Related Fields. {\bf 81}, 319-340, (1989).

\bibitem{S90}
Shiga, T.: {\it A stochastic equation based on a Poisson system
for a class of measure-valued diffusion processes}. J. Math. Kyoto
Univ. {\bf 30} (1990), 245-279.

\bibitem{SW73}
Shiga, T. and Watanabe, S., Bessel diffusions as a one-parameter
family of diffusion processes, Z. Wahrsch. verw. Geb. {\bf 27}
(1973), 37-46.

\bibitem{W86}
Walsh, J.B.: {\it An Introduction to Stochastic Partial
Differential Equations}. In: Lect. Notes Math. {\bf 1180},
265-439, Springer-Verlag (1986).

\bibitem{W97}
Wang, H.: {\it State classification for a class of measure-valued
branching diffusions in a Brownian medium}. Probab. Theory Related
Fields {\bf 109} (1997), 39-55.

\bibitem{W98}
Wang, H.: {\it A class of measure-valued branching diffusions in a
random medium}. Stochastic Anal. Appl. {\bf 16} (1998), 753-786.


\end{thebibliography}
\end{document}